\documentclass[10pt,a4wide]{article}
\pdfoutput=1

\usepackage{multicol}
\usepackage{exscale,times}
\usepackage{moredefs}

\usepackage{amssymb,amsmath,amsthm}
\usepackage{stmaryrd,mathabx}

\usepackage{tikz}
\usetikzlibrary{plotmarks}
\usetikzlibrary{shapes,arrows}
\usetikzlibrary{positioning}
\usetikzlibrary{trees,mindmap}
\usetikzlibrary{fadings,shapes.arrows,shadows}
\usetikzlibrary{decorations.pathreplacing}
\usetikzlibrary{calc}
\usetikzlibrary{patterns}

\usepackage{pgfplots}
\pgfplotsset{compat=1.10}
\usepgfplotslibrary{fillbetween}

\usepackage{graphicx}

\usepackage{lipsum}

\usepackage{hyperref}

\newcommand{\R}{\mathbb{R}}
\renewcommand{\P}{\mathbb{P}}

\usepackage{calc}
\newcommand{\WidestEntry}{$\boldsymbol C^-_{r-1}$}%
\newcommand{\SetToWidest}[1]{\makebox[\widthof{\WidestEntry}]{$#1$}}%

\begin{document}
\pagestyle{empty}

\vspace*{-3cm}
\begin{minipage}[t][1cm][t]{17.8cm}
\begin{flushright}
\fontsize{8}{10}\selectfont
\textbf{Space-Time Methods for Time-dependent Partial Differential Equations}\\
Workshop ID 1711, Mathematisches Forschungsinstitut Oberwolfach\\
Oberwolfach, March 12 - 18th, 2017
\end{flushright}
\end{minipage}

\begin{center}
%
%
{\fontsize{14}{20}\bf
Space-Time-Parallel Poroelasticity Simulation}
\end{center}

\begin{center}
\textbf{Uwe K\"ocher} $\cdot$ koecher@hsu-hamburg.de\\
\bigskip
Helmut Schmidt University, University of the Federal Armed Forces Hamburg\\
Department of Mechanical Engineering $\cdot$ Numerical Mathematics\\
Holstenhofweg 85, 22043 Hamburg, Germany\\
\end{center}
%
%



\begin{multicols}{2}
\subsubsection{Abstract}
The accurate, reliable and efficient numerical approximation of
multi-physics processes in heterogeneous porous media with varying media
coefficients
that include fluid flow and structure
interactions is of fundamental importance in
energy, environmental, petroleum and biomedical
engineering applications fields for instance.
Important applications include subsurface compaction drive, carbon sequestration,
hydraulic and thermal fracturing and oil recovery.
Biomedical applications include the simulation of
vibration therapy for osteoporosis processes of trabeculae bones,
estimating stress levels induced by tumour growth within the brain or
next-generation spinal disc prostheses.

Variational space-time methods offers some appreciable advantages such as
the flexibility of the triangulation for complex geometries in space and
natural local time stepping,
the straightforward construction of higher-order approximations and
the application of efficient goal-oriented (duality-based) adaptivity concepts.
In addition to that, uniform space-time variational methods appear to be
advantageous for stability and a priori error analyses of the discrete schemes.
Especially (high-order) discontinuous in time approaches appear to have
favourable properties due to the weak application of the initial conditions.

The development of monolithic multi-physics schemes,
instead of iterative coupling methods between the physical problems,
is a key component of the research to reduce the modeling error.
Special emphasis is on the development of efficient multi-physics and
multigrid preconditioning technologies and their implementation.

The simulation software \textsf{\footnotesize{DTM\raisebox{0.2ex}{++}}}
is a modularised framework written in \textsf{\footnotesize{C\raisebox{0.2ex}{++}11}}
and builds on top of \textsf{\footnotesize{deal.II}} toolchains.
The implementation allows parallel simulations from notebooks up to cluster scale,
cf. \cite{Koecher2017, Koecher2016a, Koecher2015}.

\vfill
\hspace*{-.15cm}
\begin{minipage}{\linewidth}
\resizebox{\linewidth}{!}{%
\footnotesize{%
\begin{tabular}{c@{\quad}c@{\quad}c}
\textsf{DTM\raisebox{0.2ex}{++}} &
\textsf{DTM\raisebox{0.2ex}{++}} &
\textsf{DTM\raisebox{0.2ex}{++}}\\
\includegraphics[height=2.3cm]{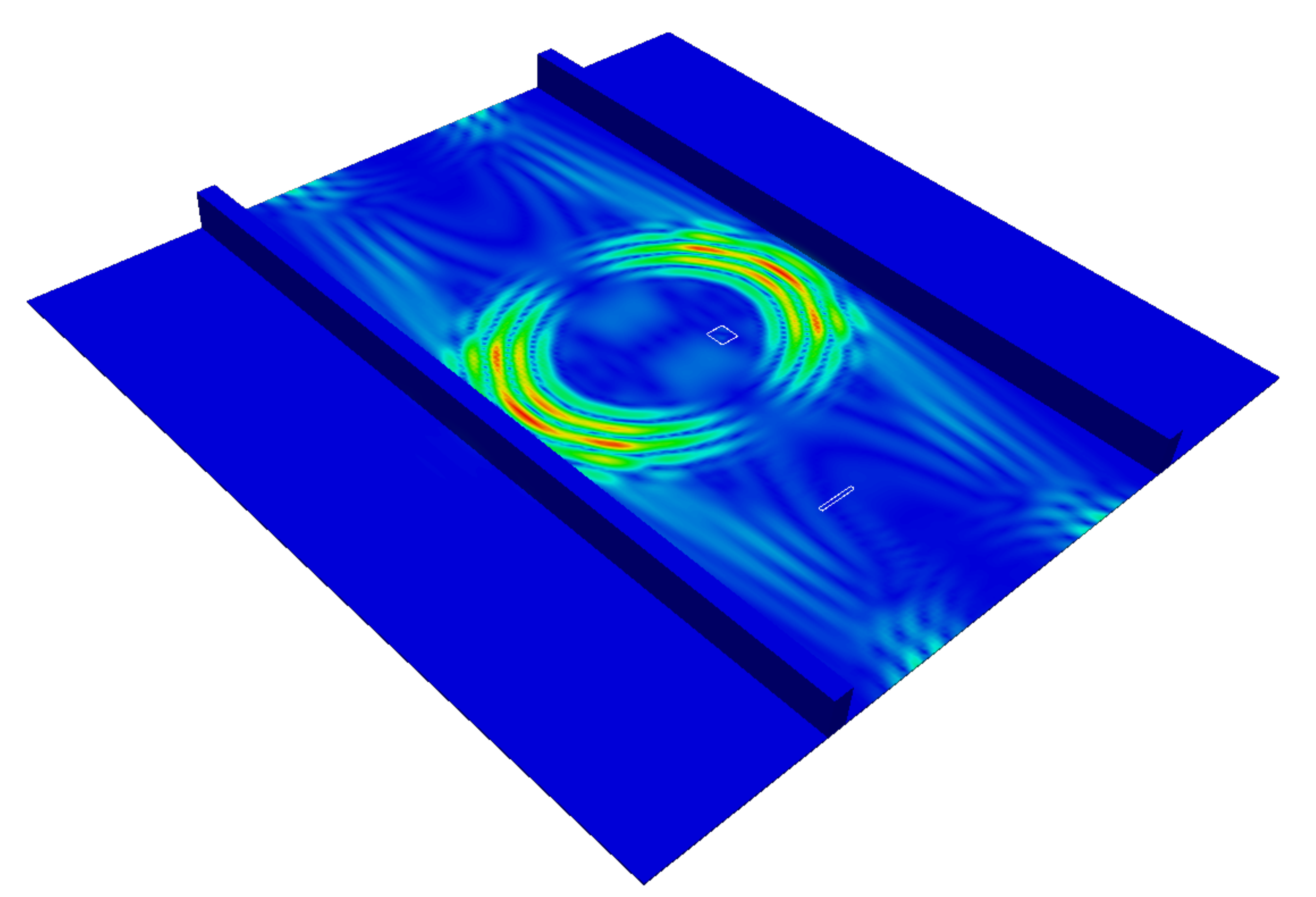} &
\includegraphics[height=2.5cm]{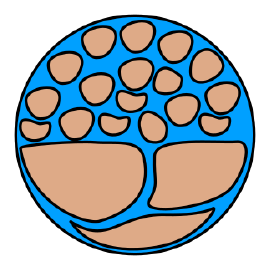} &
\includegraphics[height=2.5cm]{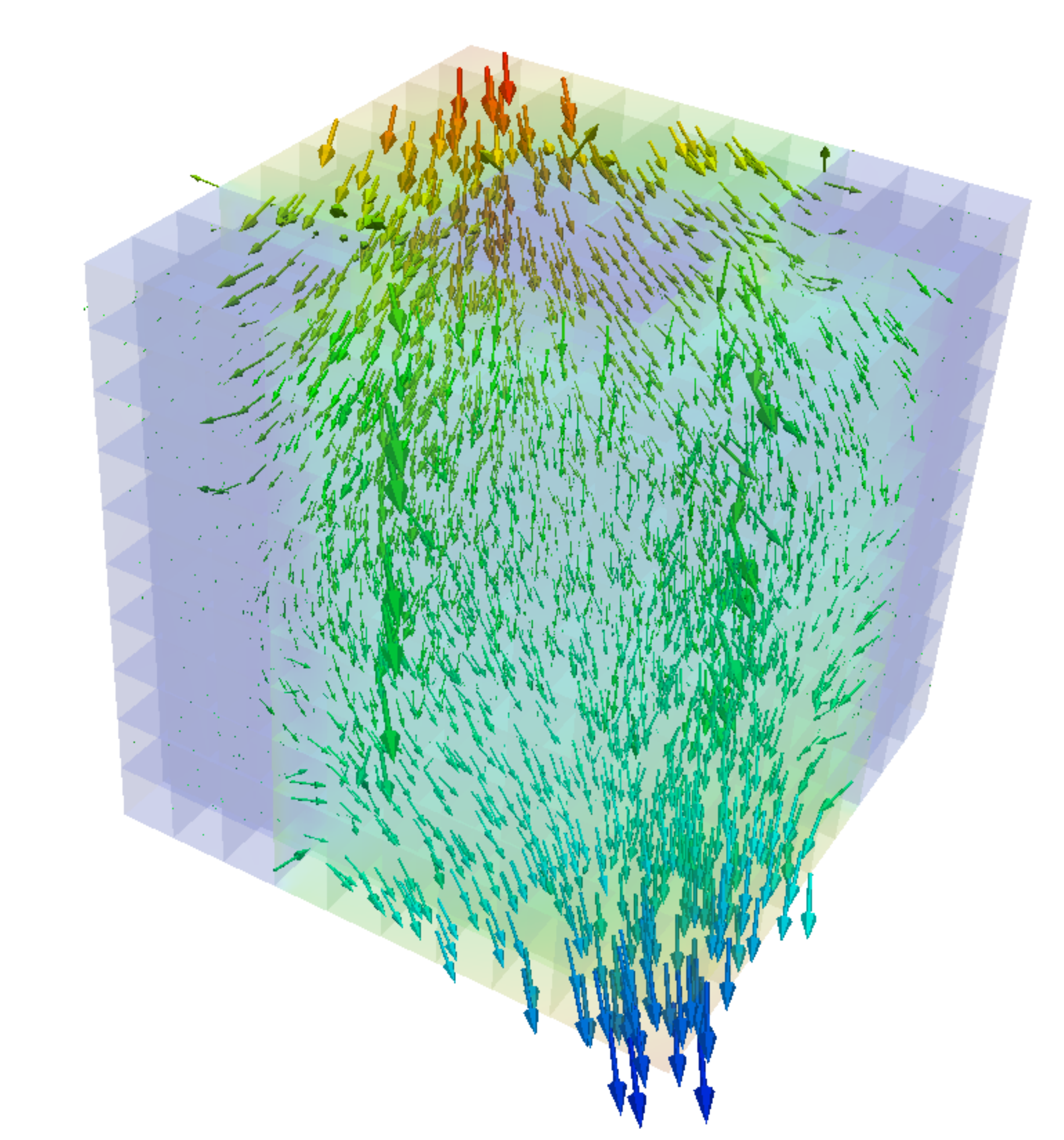}\\
\textsf{xwave} &
\textsf{biot} &
\textsf{meat}\\[1ex]
\textsf{Acoustic and Elastic Waves} &
\textsf{~} &
\textsf{Darcy-Diffusion System}\\
\end{tabular}}}
\end{minipage}
\vskip-1ex
~

\subsubsection{Poroelasticity}

The Biot poroelasticity equation system characterise a multi-physics problem
of slightly compressible single-phase diffusive porous media flow coupled with
quasi-static deformation as structure interaction.
The fully-coupled partial differential equations system in strong form for
the approximation of displacement $\boldsymbol u$, fluid flux $\boldsymbol q$
and fluid pressure $p$ is given by
\begin{displaymath}
\begin{array}{r@{\;}c@{\;}l}
-\nabla \cdot (
\boldsymbol \sigma_0
+ \boldsymbol C : \boldsymbol \epsilon(\boldsymbol u - \boldsymbol u_0)
- b\, (p - p_0) \boldsymbol 1 ) &=& \rho_b\, \boldsymbol g\,,\\[1.5ex]
\eta \boldsymbol K^{-1} \boldsymbol q + \nabla p &=& \rho_f\, \boldsymbol g\,,\\[1.5ex]
\frac{b}{K_{\textnormal{dr}}} \partial_t \sigma_v(\boldsymbol u, p) +
\big( \frac{1}{M} + \frac{b^2}{K_{\textnormal{dr}}} \big) \partial_t p +
\operatorname{div} \boldsymbol q &=& f\,,
\end{array}
\end{displaymath}
using the volumetric mean stress
\begin{displaymath}
\sigma_v(\boldsymbol u, p) = \sigma_{v,0}
+ K_{\textnormal{dr}}\, \epsilon_v( \boldsymbol u - \boldsymbol u_0 )
- b\, (p - p_0)\,,
\end{displaymath}
in $\Omega \times I$, $\Omega \subset \R^d$, $I=(0,T)$ and equipped with
appropriate initial and boundary conditions.
The coupling between the deformation and fluid flow is described by Biot using
a coupling coefficient $b = 1 - K_{\textnormal{dr}}/K_{\textnormal{s}}$, $0 < b < 1$,
and modulus $M > 0$.

\subsubsection{Discretisation in Time}

A variational discontinuous in time discretisation is applied
in order to approximate the solution by
employing Bochner spaces $\mathcal{V}$
with values in Hilbert spaces $\mathcal{H}$, e.g.
\begin{displaymath}
\begin{array}{lcl}
\mathcal{V}_{\tau}^r{(I, \mathcal{H})} &\!\!\!\!\!=&\!\!\!\!\!\Big\{
v \in L^2(I, \mathcal{H}) \,\Big|\,
v|_{I_n} \in \P_{r}(I_n, \mathcal{H})\Big\},\\[2ex]
\P_r(I_n,\mathcal{H}) &\!\!\!\!\!=&\!\!\!\!\!\Big\{
p : I_n \!\to\! \mathcal{H} \,\Big|\,
p \!= \!\sum_{j=0}^{r}{p_n^j\, t^j},\, p_n^j \in \mathcal{H}\Big\},
\end{array}
\end{displaymath}
with $I_n \!=\! (t_{n-1}, t_n)$ of a partition of the time domain.
This construction allows to solve the discrete schemes in the sense of
time marching schemes, due to the application of discontinuous Bochner spaces
as test spaces.

\subsubsection{Discretisation in Space}

The discretisation in space yields a stable, high-order,
locally mass conservative and displacement locking-free discrete scheme.
Precisely, a mixed finite element $\{\text{Nc}^f$($p$)-dQ($p$-$1$)$\}$ for
the approximation of the flux and pressure and
an interior penalty dG($p$) method
for the approximation of the displacement $\boldsymbol u$ is applied.

\subsubsection{Iterative Coupling Approach}

The Biot poroelasticity system is commonly solved with iterative coupling
approaches by the application of operator splitting methods
to decouple the fluid and structure subproblems.
A broadly accepted, unconditionally stable and fast-convergent type of such
an operator splitting method is the so called fixed-stress member.
Hereby, the flow problem is solved firstly with constrained
volumetric mean stress $\sigma_v(\boldsymbol u, p)$ and
the deformation problem is solved secondly with constrained pressure
in each iteration.
In \cite{Koecher2016a,Koecher2016} we show that an optimised fixed-stress
iterative coupling method is convergent for high-order variational time
discretisations.

\subsubsection{Monolithic Approach}

State of the art monolithic schemes commonly deploy (low-order) distributional
time integration. Arising block systems are solved with Schur-complement
technologies using standard preconditioning strategies.

Monolithic high-order dG($r$) variational time discretisations allow and need
the development of sophisticated solver and preconditioning technologies.
Spectral decomposition and an appropriate ordering of the equations yields
sparse block systems as outlined by
$$
\left[\!\!\!\!\!\!
\begin{array}{l@{\,\,}l@{\,\,}l@{\,\,}l}
\SetToWidest{\boldsymbol L_0}   & \SetToWidest{\boldsymbol C^+_1} \\[1ex]
\SetToWidest{\boldsymbol C^-_0} & \SetToWidest{\boldsymbol L_1} &
  \SetToWidest{\boldsymbol C^+_2} \\[1ex]
~ & \SetToWidest{\ddots} & \SetToWidest{\ddots} & \SetToWidest{\boldsymbol C^+_r} \\[1ex]
~ & ~ & \SetToWidest{\boldsymbol C^-_{r-1}} & \SetToWidest{\boldsymbol L_r} \\[1ex]
\end{array}
\!\!\!\!\!\!\right]
\left[\!\!\!
\begin{array}{c}
( \boldsymbol u^0, \boldsymbol q^0, p^0 )^T\\[1ex]
( \boldsymbol u^1, \boldsymbol q^1, p^1 )^T\\[1ex]
\vdots \\[1ex]
( \boldsymbol u^r, \boldsymbol q^r, p^r )^T\\[1ex]
\end{array}
\!\!\!\right]
=
\left[\!\!\!
\begin{array}{c}
\boldsymbol b^0\\[1ex]
\boldsymbol b^1\\[1ex]
\vdots \\[1ex]
\boldsymbol b^r\\[1ex]
\end{array}
\!\!\!\right]
$$
on each $I_n$ in a time marching scheme.
The coupling blocks $\boldsymbol C$ only consist of sparse contributions
from variables with time derivatives.
The diagonal block systems corresponds to lowest-order dG($0$) systems
of type
$$
\boldsymbol L_{j} =
\left[\!\!
\begin{array}{r@{\,}l@{\,\,\,\,} r@{\,}l@{\,\,\,\,} r@{\,}l}
\tau_n & \boldsymbol A & ~ & ~~\boldsymbol 0 & -b\, \tau_n & ~\boldsymbol E \\[2ex]
~& ~\boldsymbol 0 & \tau_n & \boldsymbol M_{\boldsymbol q} &
  \tau_n & ~\boldsymbol B \\[2ex]
-\lambda_{j,j}\, b & \boldsymbol E^T & \tau_n & ~\boldsymbol B^T &
  \lambda_{j,j}\, \frac{1}{M} & \boldsymbol M_p
\end{array}
\!\!\right]\,.
$$
Notably, only dG($0$) or dG($1$) type block systems must be solved,
due to the spectral decomposition, cf. \cite{Koecher2017}.
This technology facilitates parallel-in-time discrete algorithms and allows
the re-use of optimised solvers.

\subsubsection{Preconditioning technology}

An efficient preconditioning technology of the outer iterative solver GMRES for
fully-coupled dG($0$), or even dG($1$), block systems is based on an optimised
truncated fixed-stress solver technology.
Within that fixed-stress preconditioning solver, optimised solver and
preconditioning for the subproblems is used.


\subsubsection{Numerical Simulations}

\begin{minipage}{\linewidth}
\centering
\begin{tabular}{c@{\quad}c@{\quad}c}
\begin{minipage}{.25\linewidth}
\includegraphics[width=\linewidth]{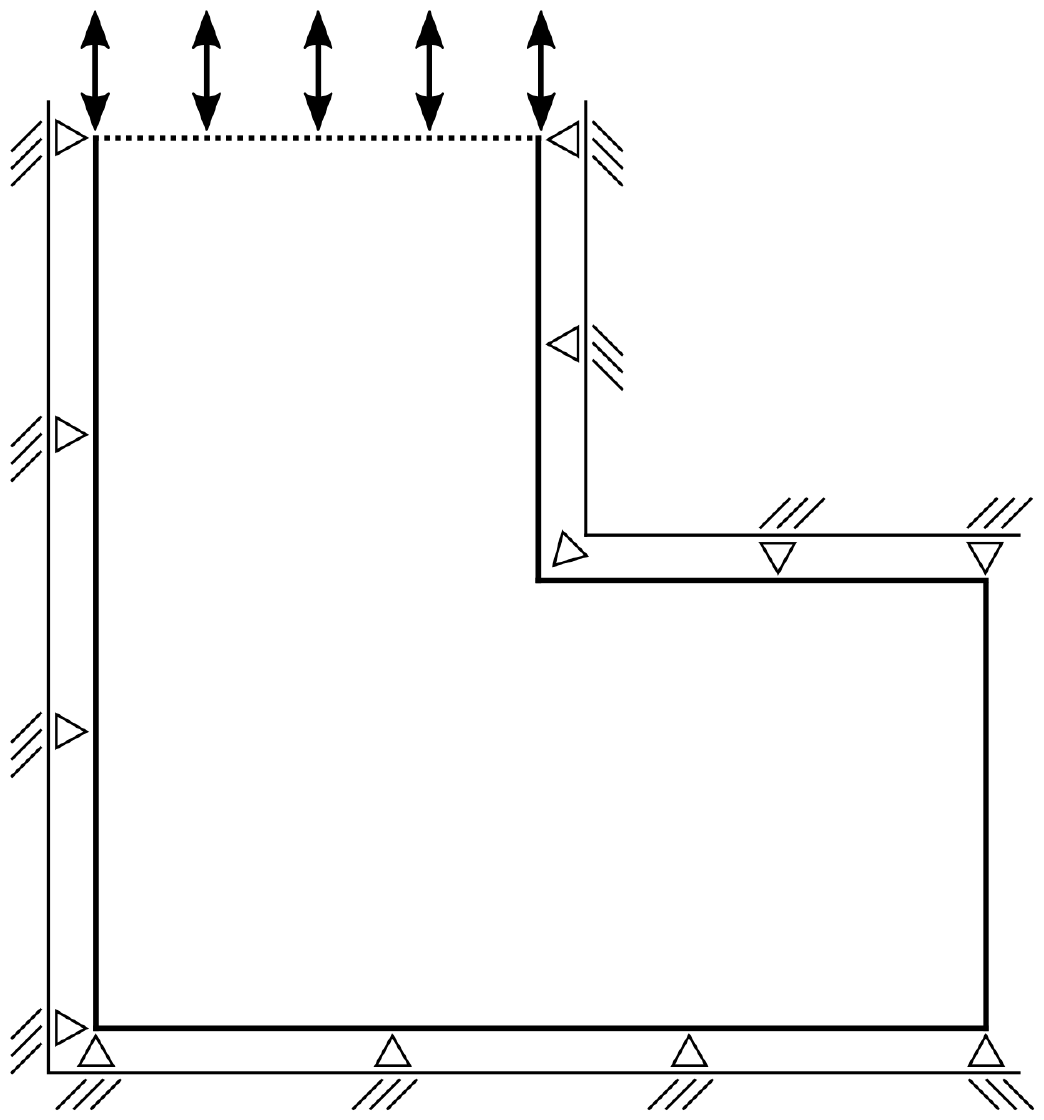}
\end{minipage} &
\begin{minipage}{.25\linewidth}
\centering
\includegraphics[width=\linewidth]{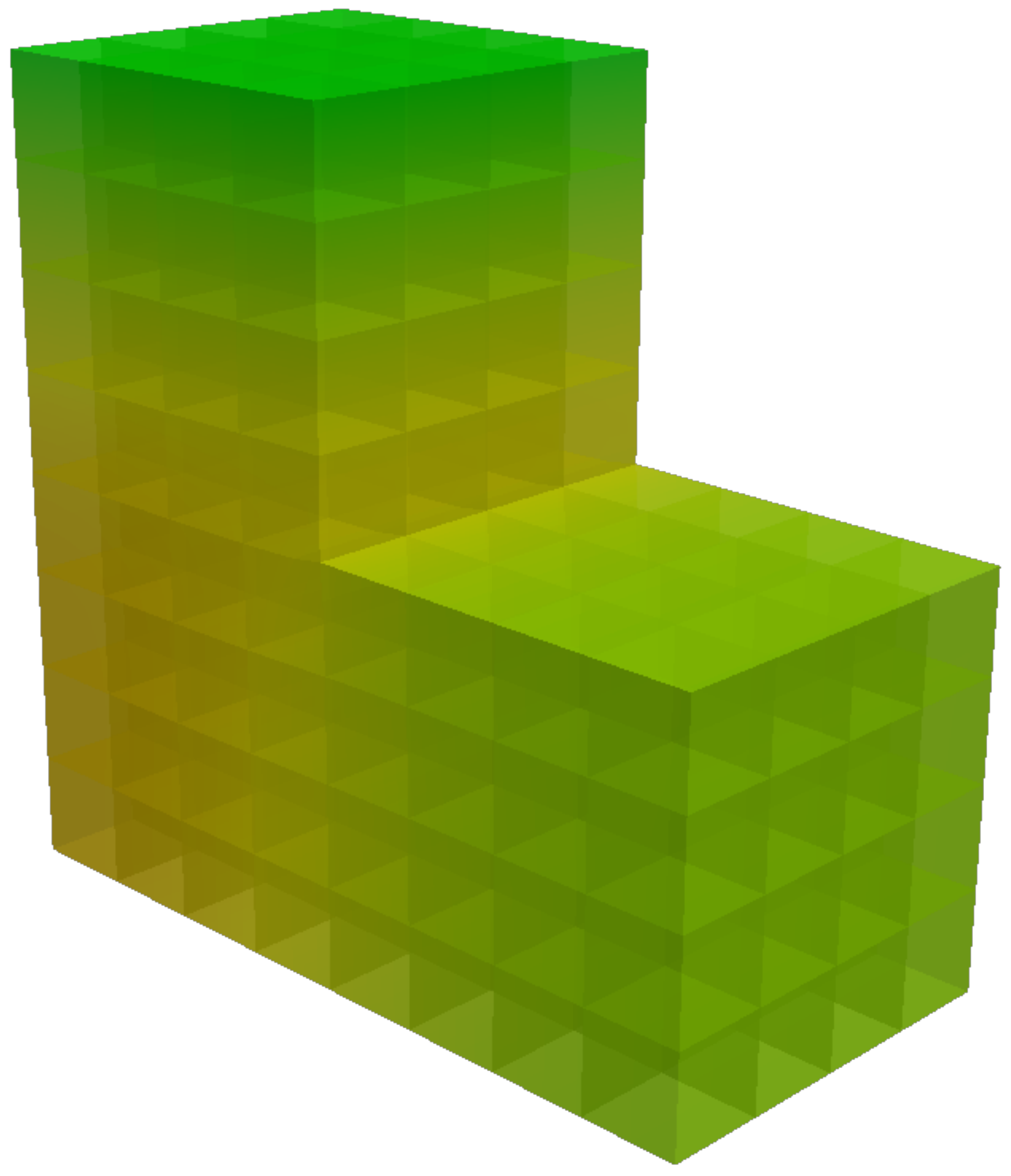}\\
\end{minipage} &
\begin{minipage}{.25\linewidth}
\centering
\includegraphics[width=\linewidth]{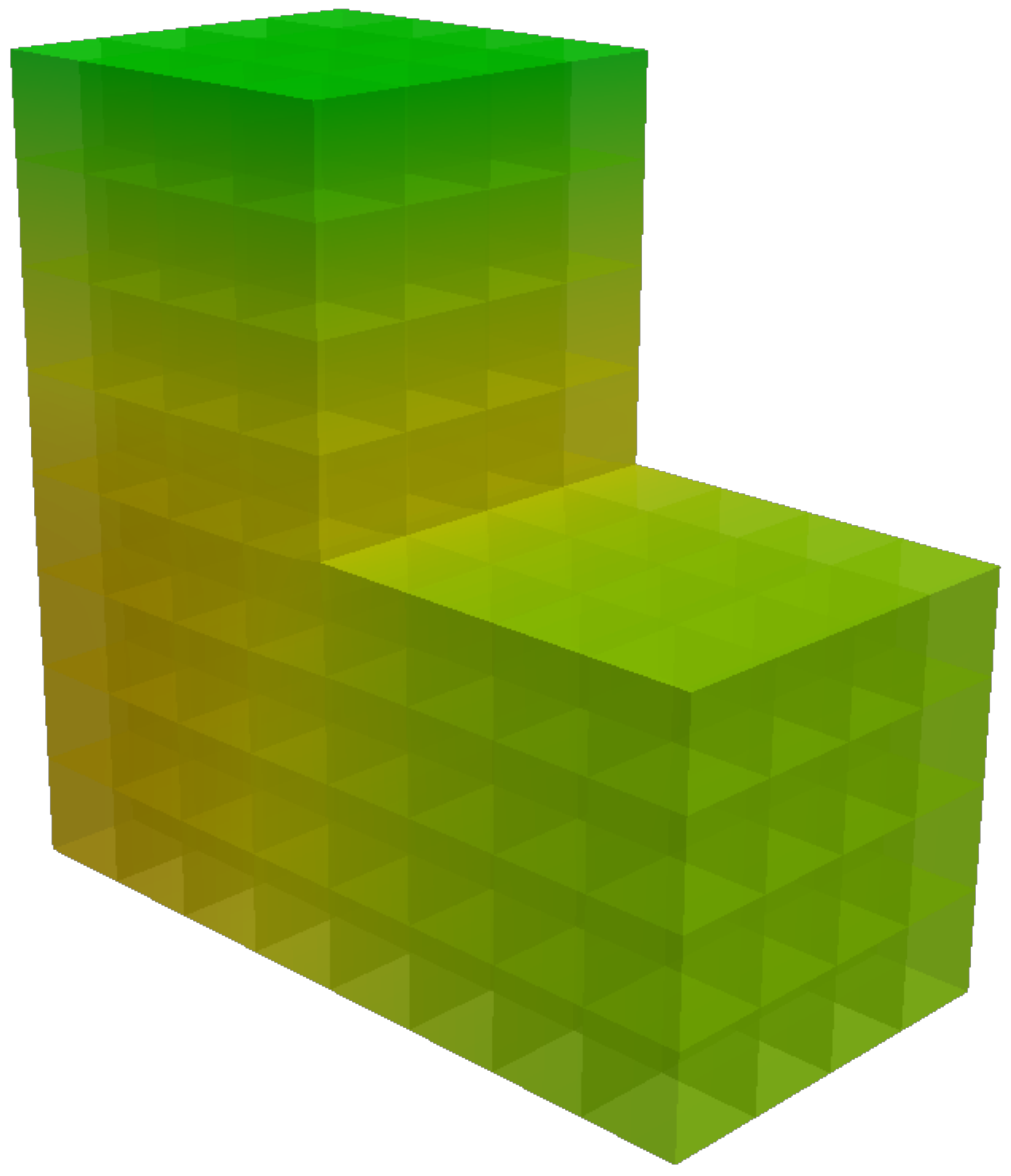}\\
\end{minipage} \\
~ & \footnotesize{dG($0$)} & \footnotesize{dG($1$)}
\end{tabular}
\end{minipage}

\begin{minipage}{\linewidth}
\centering
\footnotesize
\begin{tikzpicture}

\begin{axis}[%
width=.9\linewidth,
height=1in,
scale only axis,
/pgf/number format/.cd, 1000 sep={},
xmin=0.00,
xmax=0.60,
ymin=0,
ymax=19,
yminorticks=true,
legend style={legend pos=north east},
legend style={draw=black,fill=white,legend cell align=left},
legend entries = {
  {dG(1), $p=3$, $\tau_n=0.04$},
  {dG(0), $p=3$, $\tau_n=0.02$},
  { },
}
]

\addplot [
color=black,
solid,
line width=1.0pt,
mark=x,
mark size = 2.5,
mark options={solid,blue}
]
table[row sep=crcr]{
0.00  8\\
0.04  8\\
0.08 10\\
0.12 11\\
0.16 11\\
0.20 11\\
0.24 11\\
0.28 10\\
0.32 10\\
0.36 11\\
0.40 11\\
0.44 11\\
0.48 11\\
0.52 11\\
0.56 10\\
0.60  7\\
};

\addplot [
color=black,
solid,
line width=1.0pt,
mark=o,
mark size = 1.5,
mark options={solid,red}
]
table[row sep=crcr]{
0.00  2\\
0.02  2\\
0.04  2\\
0.06  2\\
0.08  2\\
0.10  2\\
0.12  2\\
0.14  2\\
0.16  2\\
0.18  2\\
0.20  2\\
0.22  2\\
0.24  2\\
0.26  2\\
0.28  2\\
0.30  2\\
0.32  2\\
0.34  2\\
0.36  2\\
0.38  2\\
0.40  2\\
0.42  2\\
0.44  2\\
0.46  2\\
0.48  2\\
0.50  2\\
0.52  2\\
0.54  2\\
0.56  2\\
0.58  2\\
0.60  2\\
};

\end{axis}
\end{tikzpicture}
\end{minipage}

\textbf{Fig.:} Test setting, pressure variable visualisation ($t=0.09$) and
GMRES iterations on $I=(0,T)$ for fully-coupled dG($0$) and dG($1$)
time discretisation solvers with
diagonal blockwise ($\boldsymbol L_j$) truncated fixed-stress iterative preconditioning.

\subsubsection{Future Prospects}

A challenging future development is to study the Biot-Allard model,
i.e. an extended model with dynamic structure interactions,
that includes elastic wave propagation in the structure as well as additional
memory terms.
Ongoing work is on the preconditioning technology, that is precisely the
application of multigrid techniques and cheaper discretisations,
and on their space-time-parallel implementation.

%
%

\end{multicols}

\end{document}